\documentclass{amsart}
\usepackage{amsthm,amsmath,amsfonts,amssymb,amscd,mathrsfs,graphics}

\usepackage{latexsym}
\usepackage{graphicx}

\usepackage{caption}

\usepackage{epic}

\usepackage{enumerate}
\usepackage[...]{youngtab}

\usepackage{placeins}
\usepackage{tikz}
\usepackage{listings}

\newtheorem{thm}{Theorem}[section]

\newtheorem{lem}[thm]{Lemma}
\newtheorem{prop}[thm]{Proposition}
\newtheorem{cor}[thm]{Corollary}

\makeatletter
\renewcommand{\@seccntformat}[1]{\S{\csname
the#1\endcsname}\hspace{0.5em}}
\makeatother

\begin{document}

\title{Strong Gelfand Pairs of the symplectic group  $\boldsymbol{\rm Sp}_4(q$) where $q$ is even}

\author{Stephen P.  Humphries, Joseph E. Marrow}
  \address{Department of Mathematics,  Brigham Young University, Provo, 
UT 84602, U.S.A.
E-mail:  steve@mathematics.byu.edu, jemarrow@math.byu.edu}
\date{}
\maketitle

\begin{abstract}  
A strong Gelfand pair $(G,H)$ is a finite group $G$ together with a subgroup $H$ such that every irreducible character of $H$ induces to  a multiplicity-free character of $G$.  We classify the strong Gelfand pairs of the symplectic  groups
${\rm Sp}_4(q)$ for even $q$.

\medskip

\medskip
\end{abstract}

\theoremstyle{plain}

\theoremstyle{definition}
\newtheorem*{dfn}{Definition}
\newtheorem{exa}[thm]{Example}
\newtheorem{rem}[thm]{Remark}

\newcommand{\ds}{\displaystyle}
\newcommand{\bs}{\boldsymbol}
\newcommand{\mb}{\mathbb}
\newcommand{\mc}{\mathcal}
\newcommand{\mf}{\mathfrak}
\renewcommand{\mod}{\operatorname{mod}}
\newcommand{\mult}{\operatorname{Mult}}

\def \a{\alpha} \def \b{\beta} \def \d{\delta} \def \e{\varepsilon} \def \g{\gamma} \def \k{\kappa} \def \l{\lambda} \def \s{\sigma} \def \t{\theta} \def \z{\zeta}

\numberwithin{equation}{section}

\setlength{\leftmargini}{1.em} \setlength{\leftmarginii}{1.em}
\renewcommand{\labelenumi}{\setlength{\labelwidth}{\leftmargin}
   \addtolength{\labelwidth}{-\labelsep}
   \hbox to \labelwidth{\theenumi.\hfill}}

\maketitle

\section{Introduction}

For a finite group $G$ we let $\hat G$ denote the set of irreducible characters of $G$. Then  a {\it multiplicity-free character of $G$} is a character $\chi$ of $G$ such that for   $\psi\in \hat G$, we have $\langle \chi, \psi \rangle \leq 1$. Here only complex characters are considered.
 
A {\it Gelfand pair} $(G,H)$ is a finite group $G$ together with a subgroup $H$ such that the trivial character of $H$ induces a multiplicity-free character of
$G$. The importance of Gelfand pairs is indicated by six equivalent conditions; see     \cite{AHN, BH, Har, Wik}.

A {\it strong Gelfand pair} $(G,H)$ is a finite group $G$ and   $H\le G$ such that for every   $\psi \in \hat H$  the induced character $\psi\uparrow G$ is multiplicity free. 
We will call $H$ a {\it strong Gelfand subgroup of $G$} in this situation. 
Equivalently, $(G,H)$ is a strong Gelfand pair if and only if the Schur ring determined by the {\it $H$-classes} $g^H = \{ g^h : h \in H\}, g\in G$, is commutative \cite{Har, JKar, Den}. Here  our convention is: $g^h=h^{-1}gh$. Note that $(G,G)$ is always a strong Gelfand pair.

In this paper we continue our investigation of strong Gelfand pairs of  groups that are close to being simple; in \cite {AHN, BH, GH}  we found all such pairs for 
 $G={\rm SL}(2,p^n),n\ge 1,$  $p$ a prime, and  the  symmetric groups. We refer to \cite {BH,Har, GH} for necessary background and to \cite {BC} for some of  the latest results on strong Gelfand pairs.

 We note that Gelfand pairs and strong Gelfand pairs have applications in representation theory; see  \cite {A,AA,BC,C} among many other references.
As explained above, an equivalent condition for  $(G,H)$ to be  a strong Gelfand pair is that the Schur ring determined by the $H$-classes is commutative.
This shows that a strong Gelfand pair determines a commutative Schur ring and so a commutative  association scheme, which then gives  indicates connections with algebraic combinatorics.
One other application of strong Gelfand pairs is to random walks on finite groups: if $(G,H)$ is  a strong Gelfand pair, then one can define a random walk on $G$ using probabilities that are constant on the above mentioned $H$-classes. The commutativity property of the $H$-classes means that the random walk is `diagonalizable' and so can be very well 
understood.

This paper will consider strong Gelfand pairs for the symplectic groups {\rm Sp}$_4(2^n)$ as their irreducible characters are known \cite{speven}. In contrast, the irreducible characters of ${\rm Sp}_{2k}(q),k>2,$ are not understood, where $q$ is a prime power. The  groups {\rm Sp}$_4(2^n), n>1,$ are simple and  the representation theory of these groups is considered in \cite {Da}.
 The main result of this paper is: 

\begin{thm}\label{Thm1}

The only strong Gelfand pair {\rm (Sp}$_4(2^n),H),n\ge 2,$ is where $H={\rm Sp}_4(2^n)$. 
\end{thm}

Throughout we will use the standard  Atlas notation \cite {Atl}.
\medskip

\section{Preliminary results}
\FloatBarrier

All groups considered in this paper will be assumed finite. 
For a  group $G$, the \textit{total character} of $G$, denoted $\tau_G$, is the sum of all the irreducible characters of $G$; see \cite {Tot1,Tot3,Tot2}. The following gives the `total character argument' for showing that certain subgroups are not strong Gelfand subgroups:

\begin{lem} [Lemma $3.3$ \cite {GH}] \label{lem:totalchar} Let $H\leq G$ be   groups. If there is  $\chi\in \hat G$ with 
$$
\deg(\tau_H) < \deg(\chi),
$$
 then $(G, H)$ is not a strong Gelfand pair.
\end{lem}

The following indicates that it is important to determine which maximal subgroups are strong Gelfand pairs.

\begin{lem}  [Lemma $3.1$ \cite {BH}]  \label{lem:totalstack} Suppose we have groups $H\leq K\leq G$. If $(G, K)$ is not a strong Gelfand pair, then neither is $(G, H)$.
\end{lem}

For  $q=2^e$,  $e>1$, we find from Table $8.14$ of \cite{max} that   the maximal subgroups of $\mathrm{Sp}_4(q)$  are as listed in Table \ref{maxeven}.

\begin{table}[h!]
\centering
\caption{Maximal subgroups of $\mathrm{Sp}_4(q)$, for $q=2^e, e>1$}
\label{maxeven}
\vspace{5pt}
\begin{tabular}{|ccc|}
\hline
Group & Order& Conditions\\
\hline
$\vphantom{E^{3^2}}E_q^3\colon\mathrm{GL}_2(q)$ & $q^3\cdot (q^2+q)(q-1)^2$ & \\%
$\vphantom{E^{3^2}}E_q^3\colon\mathrm{GL}_2(q)$ & $q^3\cdot (q^2+q)(q-1)^2$ & \\%
$\mathrm{Sp}_2(q)\wr 2$ & $2q^2(q^2-1)^2$ & \\%
$\mathrm{Sp}_2(q^2)\colon 2$ & $2q^2(q^4-1)$ & \\
$\mathrm{Sp}_4(q_0)$ & $q_0^4(q_0^2-1)(q_0^4-1)$ & $q=q_0^r, r$ is prime\\%
$\mathrm{SO}_4^+(q)$ & $2q^2(q^2-1)^2$ & \\
$\mathrm{SO}_4^-(q)$ & $2q^2(q^4-1)$ & \\
$\mathrm{Sz}(q)$ & $q^2(q^2+1)(q-1)$ & $e$ odd\\%
\hline
\end{tabular}
\end{table}

Table \ref{maxeven} has the maximal  subgroup $E_q^3 \colon \mathrm{GL}_2(q)$ listed twice because there are two  non-conjugate maximal subgroups of Sp$_4(q)$ which are isomorphic to $E_q^3 \colon \mathrm{GL}_2(q)$.

From Lemma \ref {lem:totalstack} Theorem \ref {Thm1} will follow if we can show that none of these maximal subgroups is a strong Gelfand subgroup. We consider each case separately.

The next two results  will allow us to assume   $e\ge 3$.   

We first  consider  the symplectic group ${\rm Sp}_4(2)$; since ${\rm Sp}_4(2)\cong S_6$ the result here follows from our consideration of the symmetric groups in 
\cite {AHN}:

\begin{prop} \label{thm:sp2}
The only proper subgroups of $\mathrm{Sp}_4(2)$ which are strong Gelfand subgroups are  the maximal subgroups.\qed
\end{prop}

\begin{prop}\label{cor:sp4}
No proper subgroup of $\mathrm{Sp}_4(4)$ is a strong Gelfand subgroup.
\end{prop}
\noindent {\it Proof} 
 We use the MAGMA \cite {Mag} code given in the Appendix  to obtain  this result.\qed\medskip
 
 In what follows we will often have the situation where $H \le G, |G:H|=2$. We introduce the following conventions.
 For $\psi \in \hat H$ it  is well-known \cite {JaL} that either 
 
 \noindent (i) $\psi \uparrow G$ is a sum of two distinct characters  in $\hat G$  (call this the {\it splitting case}); or
 
 \noindent (ii)  $\psi \uparrow G$ is irreducible (call this the {\it fusion} case). 
 
 In the splitting case, if 
  $\psi \uparrow G=\chi_1+\chi_2, \chi_1,\chi_2 \in \hat G$, then $\chi_i\downarrow H=\psi, i=1,2.$
 
\medskip

In the situation $|G:H|=2$ the relationship between $\tau_G$ and $\tau_H$ is given in: 
\begin{lem}\label{lemsplitf}
Let $H\le G, |G:H|=2$. Let $\mathcal S$ be the set of $\psi \in\hat H$ that split and 
let  $\mathcal F$ be the set of $\psi \in\hat H$ that fuse. Then 
$$
\tau_G(1)=2\sum_{\psi \in \mathcal S}\psi(1)+  \sum_{\psi \in \mathcal F}\psi(1).\qed 
$$
\end{lem}

The character table for $\mathrm{Sp}_4(q)$ is given in \cite{speven} and  we will use notation from \cite{speven}. 

\begin{thm}  \cite{speven}.\label{lem:totalsymplectic}  (i) The degree of the total character of $\mathrm{Sp}_4(q)$ is $q^6+q^4-q^2$ if $q$ is even.

(ii) The largest degree of an irreducible character of $\mathrm{Sp}_4(q)$ is $q^4+2q^3+2q^2+2q+1$ when $q\geq 4$ is a power of $2$.
\end{thm}
\noindent\textit{Proof.} (i) We just  sum  the degrees of characters of $\mathrm{Sp}_4(q)$ as listed in \cite{speven}. (ii) follows directly from \cite {speven}.\qed\medskip

\begin{lem}\label{isomorpicmaximal}
If $q=2^e, e>1$, then $\mathrm{Sp}_2(q)\wr 2 \cong \mathrm{SO}_4^+(q)$  and $\mathrm{Sp}_2(q^2)\colon 2 \cong \mathrm{SO}_4^-(q)$.
\end{lem}
\noindent\textit{Proof.} 
See Proposition $7.2.1$  and Table $8.14$ of \cite{max}.  \qed\medskip

We now consider the maximal subgroups separately in the following sections.

\section { The maximal subgroup $\mathrm{Sp}_2(q)\wr 2 $}


\begin{thm}\label{wreath}
For $q=2^e, e>1,$ the maximal subgroup $\mathrm{Sp}_2(q)\wr 2 \leq \mathrm{Sp}_4(q)$ is not a strong Gelfand subgroup.
\end{thm}
\noindent\textit{Proof.} This proof will be a `total character argument' and so we will need to find the total character of
$\mathrm{Sp}_2(q)\wr 2$. 
We have $\mathrm{Sp}_2(q)\wr 2 = \mathrm{Sp}_2(q)^2: 2$ and one way to  represent the elements of $\mathrm{Sp}_2(q)^2: 2$ is by $2 \times 2$ blocks of $2 \times 2$ matrices, where the cyclic  subgroup $2$ is generated  by  $\begin{bmatrix} 0&I_2\\I_2&0\end{bmatrix}$ and
$(a,b) \in \mathrm{Sp}_2(q)^2$ is represented as the block matrix 
$\begin{bmatrix} a&0\\0&b\end{bmatrix}.$

Now $\mathrm{Sp}_{2}(q) \cong \mathrm{SL}_2(q)$ has  character table  given in \cite{Dor}  (see also  \cite {GH}); we reproduce it here in Table \ref{sleven}.
Here the parameters $s,t,j,m$ satisfy $1\leq s, t \leq (q-2)/2$, $1\leq j, m \leq q/2$,  $\rho$ is a primitive $(q-1)$-th root of unity and $\sigma$ a primitive $(q+1)$-th root of unity. 

\begin{table}[h!]
\centering
\caption{Character Table for $\mathrm{SL}_{2}(q)$ with $q$ even}
\label{sleven}
\vspace{5pt}
\begin{tabular}{c|cccc}
Class & $1$ & $c$ & $a^t$ & $b^m$\\
Size & $1$ & $q^2-1$ & $q(q+1)$ & $q(q-1)$\\
\hline
 Tr & $1$ & $1$ & $1$ & $1$\\
$\psi$ & $q$ & $0$ & $1$ & $-1$\\
$\chi_s$ & $q+1$ & $1$ & $\rho^{st}+\rho^{-st}$ & $0$\\
$\theta_j$ & $q-1$ & $-1$ & $0$ & $-(\sigma^{jm}+\sigma^{-jm})$\\
\end{tabular}
\end{table}

Here the conjugacy classes of SL$_2(q)$ are represented by powers of the following elements:
$$
1 = \begin{bmatrix}
1 & 0\\
0 & 1
\end{bmatrix},\ 
c=\begin{bmatrix}
1 & 0\\
1 & 1
\end{bmatrix},\ 
a=\begin{bmatrix}
\rho & 0\\
0 & \rho^{-1}
\end{bmatrix},
$$
and an element $b$ of order $q+1$. We also give the sizes of the classes in Table \ref{sleven}.

Since $\mathrm{Sp}_2(q)\wr 2 \cong \mathrm{Sp}_2(q)^2\colon 2$
  the irreducible characters of $\mathrm{Sp}_2(q)\wr 2$ are easily found using Table 2.
In Table 
   3 we give the degrees of the irreducible characters of $\mathrm{Sp}_2(q)\wr 2$.
These character degrees are obtained using \cite[Proposition $20.9$, Theorem $19.18$]{JaL}. Further, in Table 
 3 we are assuming that  
 $$1 \leq s, s^\prime \leq (q-2)/2, \,\,1 \leq j, j^\prime \leq q/2 \text { and } s\neq s^\prime,j\neq j^\prime.
 $$

\begin{table}[h!]
\label{chardegwreath}
\centering
\caption{Character degrees for $\mathrm{Sp}_2(q)\wr 2$ with $q$ even}
\vspace{5pt}
\begin{tabular}{c|cc}
Character & Degree & Multiplicity\\
\hline
$\vphantom{()^{2^2}}(\mathrm{ Tr}\times\mathrm{ Tr})_1$ & $1$ & $1$\\
$(\mathrm{ Tr}\times\mathrm{ Tr})_2$ & $1$ & $1$\\
$\mathrm{ Tr}\times\psi$ & $2\cdot q$ & $1$\\
$\mathrm{ Tr}\times\chi_i$ & $2\cdot(q+1)$ & $(q-2)/2$\\
$\mathrm{ Tr}\times\theta_j$ & $2\cdot (q-1)$ & $q/2$\\
$(\psi\times\psi)_1$ & $q^2$ & $1$\\
$(\psi\times\psi)_2$ & $q^2$ & $1$\\
$\psi\times\chi_s$ & $2\cdot q(q+1)$ & $(q-2)/2$\\
$\psi\times\theta_j$ & $2\cdot q(q-1)$ & $q/2$\\
$(\chi_s\times\chi_s)_1$ & $(q+1)^2$ & $(q-2)/2$\\
$(\chi_s\times\chi_s)_2$ & $(q+1)^2$ & $(q-2)/2$\\
$\chi_s\times\chi_{s^\prime}$ & $2\cdot (q+1)^2$ & $(q-2)(q-4)/8$\\
$\chi_s\times\theta_j$ & $2\cdot(q^2-1)$ & $q(q-2)/4$\\
$(\theta_j\times\theta_j)_1$ & $(q-1)^2$ & $q/2$\\
$(\theta_j\times\theta_j)_2$ & $(q-1)^2$ & $q/2$\\
$\theta_j\times\theta_{j^\prime}$ & $2\cdot (q-1)^2$ & $q(q-2)/8$
\end{tabular}
\end{table}

In Table 3 the suffices $1,2$ are written to indicate that these  are the split cases. The lack of such a suffix indicates the fusion cases. 
In Table 3 each case has a certain `Multiplicity' that is also indicated; this  depends on the parameters involved. 
Then from 
 Table 3 we obtain the degree of the total character of  Sp$_2(q)\wr 2$:
\begin{align*}
(\tau_{\mathrm{Sp}_2(q)\wr 2})(1)
&=1 + 1 + 2 q + (q + 1) (q - 2) + (q - 1) \cdot q + 2 q^2 + q\cdot(q + 1)\cdot (q - 2)\\ & \qquad + q^2 (q - 1)
 + (q + 1)^2 \cdot(q - 2) + (q + 1)^2\cdot (q - 2)\cdot (q - 4)/4\\& \qquad + (q^2 - 1)\cdot q\cdot (q - 2)/2 + (q - 1)^2\cdot q
 + (q - 1)^2\cdot q\cdot (q - 2)/4 \\&
= q^4+q^3-q.
\end{align*}

Now $q^4+q^3-q < q^4+2q^3+2q^2+2q+1$, and  by Theorem  \ref{lem:totalsymplectic}   $q^4+2q^3+2q^2+2q+1$ is the degree of an irreducible character of $\mathrm{Sp}_4(q)$. Then by Lemma \ref{lem:totalchar}  $(\mathrm{Sp}_4(q), \mathrm{Sp}_2(q)\wr 2)$ is not a strong Gelfand pair.\qed\medskip 

By Lemma \ref {isomorpicmaximal} and the fact that the above argument is a `total character argument' (not dependent on the particular embedding of $ \mathrm{Sp}_2(q)\wr 2$ in $\mathrm{Sp}_4(q)$)        we see that we have now also dealt with  the maximal subgroup ${\rm SO}^+_4(q)$ case from Table 1:

\begin{cor} \label{cor1} The maximal subgroup {\rm SO}$^+_4(q)<{\rm Sp}_4(q)$ is not a strong Gelfand subgroup.\qed
\end{cor}\medskip

\section {The maximal subgroups $E_q^3:{\rm GL}_2(q)$}

By Theorems \ref {thm:sp2} and \ref {cor:sp4}  we may assume that $q>4$.
\begin{thm}\label{thm:elementaryabelian}
For $q=2^e, e>2$, the maximal subgroup $E_q^3\colon\mathrm{GL}_2(q) \leq \mathrm{Sp}_4(q)$ is not a strong Gelfand subgroup.
\end{thm}
\noindent\textit{Proof.} In \cite{speven} two isomorphic maximal subgroups are considered; they are denoted $P$ and $Q$. 
The orders of $P$ and $Q$ are $q^3(q^2+q)(q-1)^2$ and they are isomorphic to  $E_q^3\colon\mathrm{GL}_2(q)$. The character tables for these subgroups are given in \cite{speven}.

We  take the inner product of the character of $P$ denoted by  $\chi_5(k)$ in  \cite{speven} with a character of $\mathrm{Sp}_4(q)$ restricted to $P$, namely $\chi_{1}(m, n)\downarrow P$. In what follows  $A_i,A_{ij}, C_j, D_k$ is the notation used  in \cite{speven} for the classes of $P$; further, the sizes of these classes are also given in \cite{speven}. Using all of this information we obtain:
\begin{align*}
&\langle \chi_5(k), \chi_1(m, n)\downarrow P\rangle\\&
=\frac{1}{|P|}\Big(
|A_1| q(q^2-1) (q+1)^2(q^2+1)
+|A_2| q(q-1) (q+1)^2 
+|A_{31}| (-q)(q+1) (q+1)^2\\
&\qquad + |A_{32}| (-q) (2q+1)
+|C_2(i)| (q-1)\alpha_{ik} (q+1)\alpha_{im}\alpha_{in} +|D_2(j)| (-\alpha_{jk}) \alpha_{jm}\alpha_{jn}
\Big)\\&
=\frac{1}{q^4(q-1)(q^2-1)}\Big(q(q^2-1) (q+1)^2(q^2+1)\\
&\qquad +(q^2-1) q(q-1) (q+1)^2
+(q-1) (-q)(q+1) (q+1)^2\\ 
&\qquad +(q-1)(q^2-1) (-q) (2q+1)\\ 
&\qquad +\sum_{i=1}^{(q-2)/2} q^3(q+1) (q-1)\alpha_{ik} (q+1)\alpha_{im}\alpha_{in}
+\sum_{j=1}^{(q-2)/2} q^3(q^2-1) (-\alpha_{jk}) \alpha_{jm}\alpha_{jn}\Big)\\
&=\frac{1}{q^7-q^6-q^5+q^4}\Big(q^7+2q^6+q^5-q^3-2q^2-q
 +q^6+q^5-2q^4-2q^3+q^2\\
 &\qquad \qquad \qquad +q
-q^5-2q^4+2q^2+q-2q^5+q^4+3q^3-q^2-q\\
&\qquad +\sum_{i=1}^{(q-2)/2} (q^6+q^5-q^4-q^3)\alpha_{ik}\alpha_{im}\alpha_{in}
+\sum_{j=1}^{(q-2)/2}(-q^5+q^3)\alpha_{jk}\alpha_{jm}\alpha_{jn}\Big)\\
\end{align*}
\begin{align*}&=
\frac{1}{q^7-q^6-q^5+q^4} \left(q^7+3q^6-q^5-3q^4+\left(q^6-q^4\right)\sum_{j=1}^{(q-2)/2} \alpha_{jk}\alpha_{jm}\alpha_{jn}\right)\\&
=\frac{3+q}{q-1}+\frac{1}{q-1}\left(\sum_{j=1}^{(q-2)/2} \alpha_{jk}\alpha_{jm}\alpha_{jn}\right).
\end{align*}

Here $\alpha_{ij} = \overline{\gamma}^{ij}+\overline{\gamma}^{-ij}$ where $\langle\gamma\rangle = \mathbb{F}_q^\times$, and $\overline{\gamma}$ is the image of $\gamma$ under a fixed monomorphism from $\mathbb{F}_q^\times$ into $\mathbb{C}^\times$, making $\overline{\gamma}$ a $(q-1)$-th root of unity. For clarity of notation, in our calculations we omit the overline.

Now supposing that $q>5$, if we have $\sum_{i=1}^{(q-2)/2}\alpha_{ik}\alpha_{im}\alpha_{in}=q-5$, then the above gives $\langle \chi_5(k), \chi_1(m, n)\downarrow P\rangle=2$. We will now show that there is a choice of $k, m, n$ so that $\sum_{i=1}^{(q-2)/2}\alpha_{ik}\alpha_{im}\alpha_{in}$ is equal to  $q-5$.
We calculate:
\begin{align*}
&\sum_{j=1}^{(q-2)/2}\alpha_{jk}\alpha_{jm}\alpha_{jn}\\
&= \left(\sum_{j=1}^{(q-2)/2} \gamma^{jk}+\gamma^{-jk}\right)\left(\sum_{j=1}^{(q-2)/2} \gamma^{jm}+\gamma^{-jm}\right)\left(\sum_{j=1}^{(q-2)/2} \gamma^{jn}+\gamma^{-jn}\right)\\
&=\sum_{j=1}^{q-2} \gamma^{j(k+m+n)} + \sum_{j=1}^{q-2} \gamma^{j(k+m-n)} + \sum_{j=1}^{q-2} \gamma^{j(k-m+n)} + \sum_{j=1}^{q-2} \gamma^{j(k-m-n)}
\end{align*}
and notice that each of these four sums will be $q-2$ if $q-1$ divides  $j$, and $-1$ otherwise. Suppose that $q>5$ and choose $k=q-4, m=1, n=2$. Then  $m\neq n$ and $m+n\neq q-1$, as required. We also have that only one of $k\pm m \pm n $ is congruent to zero mod
  ${q-1}$. This gives
\begin{align*}
&\sum_{j=1}^{q-2} \gamma^{j(q-1)} + \sum_{j=1}^{q-2} \gamma^{j(q-3)} + \sum_{j=1}^{q-2} \gamma^{j(q-5)} + \sum_{j=1}^{q-2} \gamma^{j(q-7)} = q-5.
\end{align*}

Then for $k=q-4, m=1, n=2$ we have:
\begin{align*}
&\langle \chi_5(q-4), \chi_1(1, 2)\downarrow P\rangle = \frac{3+q+\left(\sum_{j=1}^{(q-2)/2} \alpha_{j(q-4)}\alpha_{j1}\alpha_{j2}\right)}{q-1}
 = \frac{3+q+q-5}{q-1}
= 2,
\end{align*}
showing that $(\mathrm{Sp}_4(q), P)$ is not a strong Gelfand pair if $q>5$.

A similar argument shows that $(\mathrm{Sp}_4(q), Q), q>5,$ is also not a strong Gelfand pair.
\qed
\medskip

\section {The maximal subgroups Sp$_2(q^2):2$ and {\rm Sp}$_4(q_0)$}

The elements of  the field $\mathbb F_{q^2}$ can be represented as $2 \times 2$ matrices over $\mathbb F_q$. This shows how Sp$_2(q^2)\le
{\rm Sp}_4(q)$. The action of the $2$ in Sp$_2(q^2):2$ is the Galois action. 

\begin{thm}
For $q=2^e, e>1,$ the pair $\left(\mathrm{Sp}_4(q), \mathrm{Sp}_2\left(q^2\right)\colon 2\right)$ is not a strong Gelfand pair.
\end{thm}
\noindent\textit{Proof.} Let $G=\mathrm{Sp}_2\left(q^2\right)\colon 2$ and $H = \mathrm{Sp}_2\left(q^2\right)\leq G$. Using  Table \ref{sleven} we get the character table for $H$; see  Table \ref{lineardegrees} where  $1 \leq s \leq (q^2-2)/2$ and $1 \leq j \leq q^2/2$.

\begin{table}[h!]
\centering
\caption{Character degrees for $\mathrm{Sp}_{2}\left(q^2\right)$ with $q$ even}
\vspace{5pt}
\label{lineardegrees}
\begin{tabular}{c|c c}
Character & Degree&Multiplicity\\
\hline
 Tr & $1$& $1$\\
$\psi$ & $q^2$& $1$\\
$\chi_s$ & $q^2+1$ & $(q^2-2)/2$\\
$\theta_j$ & $q^2-1$&$q^2/2$
\end{tabular}
\end{table}

In order to find the degree of $\tau_{G}$, we will need to determine which characters of $H$ split and which  fuse; 
it will suffice to determine which characters of $H$ induce to  irreducible characters of $G$. Again from \cite{JaL}, since $|G \colon H|=2,$ we know that, by  inducing, every character in $\hat{H}$ either splits into a sum of two irreducible characters or fuses pairwise into irreducible characters in $\hat{G}$.
We  use Lemma \ref {lemsplitf} and Table 4 to give:

\begin{prop}\label{propsplit} Let $G={\rm Sp}_2(q^2):2\ge H={\rm Sp}_2(q^2)$. Then 

\noindent (i) $\text{Tr}_H$ splits;

\noindent (ii)  $\psi$ splits;

\noindent (iii) all $\theta_j$ fuse;

 The characters $\chi_s$ sometimes split, but not always: 
 
\noindent (iv)  $\chi_s\uparrow G$ is  irreducible  if $(q^2-1)\nmid s(q\pm1)$; and

 \noindent (v) $\chi_s\uparrow G$ is the sum of two irreducible characters if $(q^2-1) \mid s(q\pm1)$.
\end{prop}
\noindent{\it Proof} (i) It is clear that Tr$_H$ splits.

\noindent (ii) Since $\psi(1)=q^2$ and there is no other character of degree $q^2$ we see that $\psi$ cannot fuse.

\noindent (iii)  It will suffice to show that $\langle \theta_j \uparrow G,\theta_j \uparrow G\rangle=1$. 
Now a calculation shows that $\theta_j \uparrow \mathrm{Sp}_2(q^2)\colon 2$ 
is as  described in the following  table, where $\sigma$ is a primitive $(q^2+1)$-th root of unity. 
\\[5mm]
\begin{tabular}{c|ccccc}
 & $\text{Tr}_H$ & $c$ & $a^t$ & $b^m$ & $G \smallsetminus H$\\
\hline\\
$\theta_j \uparrow \mathrm{Sp}_2(q^2)\colon 2$ & $2q^2-2$ & $-2$ & $0$ & $-(\sigma^{jm} + \sigma^{-jm}+\sigma^{jmq} + \sigma^{-jmq})$ & $0$
\end{tabular}
\\[5mm]
Now 
$(\theta_j\uparrow G)(G \setminus H)=\{0\}$ and 
for $g \in H$ we have $g$ and $g^{-1}$ are conjugate. Thus 
\begin{align*}
\langle \theta_j\uparrow G, \theta_j\uparrow G\rangle &= \frac{1}{|G|} \sum_{g\in G}(\theta_j\uparrow G)(g)\cdot (\theta_j\uparrow G)\left(g^{-1}\right)\\&=
\frac{1}{|G|} \sum_{g\in H}(\theta_j\uparrow G)(g)\cdot (\theta_j\uparrow G)(g^{-1})
=\frac{1}{|G|} \sum_{g\in H}(\theta_j\uparrow G)^2(g).
\end{align*}

Using Table 2 again and taking $g_m \in (b^m)^G$  the above is equal to
\begin{align*}&
\frac{1}{2(q^6-q^2)}\left(\underbrace{(2q^2-2)^2}_{\text{ Tr}} + \underbrace{(-2)^2(q^4-1)}_{c} + \underbrace{0}_{\text{$a^t$}}+\underbrace{(q^4-q^2)}_{\text{size of $(b^m)^H$}} {\sum_{m=1}^{q^2/2}}\left(\theta_j\uparrow G\right)^2(g_m)\right)
\\&=
\frac{1}{2(q^6-q^2)}\left({8q^4-8q^2} + {(q^4-q^2)} {\sum_{m=1}^{q^2/2}}  {\left(\theta_j\uparrow G\right)^2(g_m)}\right)
\\&=
\frac{1}{2(q^6-q^2)}\left(8q^4-8q^2+(q^4-q^2)\sum_{m=1}^{q^2/2}\left(-\sigma^{jm}-\sigma^{-jm}-\sigma^{jmq}-\sigma^{-jmq}\right)^2\right)
\\&=
\frac{1}{2(q^6-q^2)}
\bigg (
8q^4-8q^2+(q^4-q^2)\sum_{m=1}^{q^2/2}
\big (
4 +  {(\sigma^{2jm}+\sigma^{-2jm})}+ (\sigma^{2jmq} +\sigma^{-2jmq})\\&\qquad  + (2\sigma^{jm(q+1)} + 2\sigma^{-jm(q+1)}) + (2\sigma^{jm(q-1)} + 2\sigma^{-jm(q-1)})
\big )
\bigg )
\end{align*}

Now, since $1+\sum_{i=1}^{q^2/2} \sigma^i+\sigma^{-i}=\sum_{i=0}^{q^2} \sigma^i$,  the above  is
\begin{align*}&
\frac{1}{2(q^6-q^2)}\left(8q^4-8q^2+(q^4-q^2)\sum_{m=1}^{q^2}\left(
2+\sigma^{2jm}+\sigma^{2jmq}+2\sigma^{jm(q+1)}+2\sigma^{jm(q-1)}
\right)\right)\\&\qquad = 
\frac{(8q^4-8q^2+2q^2(q^4-q^2)-6(q^4-q^2))}{2q^6-2q^2} = \frac{2q^6-2q^2}{2q^6-2q^2} = 1
\end{align*}
as required for (iii). 
\medskip

\noindent (iv) 
Now a calculation shows that $\chi_s \uparrow \mathrm{Sp}_2(q^2)\colon 2$ 
is as  described in the following  table, where $\rho$ is a primitive $(q^2-1)$-th root of unity. 
\\[5mm]
\begin{tabular}{c|ccccc}
 & $\text{Tr}_H$ & $c$ & $a^t$ & $b^m$ & $G \smallsetminus H$\\
\hline\\
$\theta_j \uparrow \mathrm{Sp}_2(q^2)\colon 2$ & $2q^2+2$ & $2$ &  $-(\rho^{jm} + \rho^{-jm}+\rho^{jmq} + \rho^{-jmq})$ &$0$  & $0$
\end{tabular}
\\[5mm]
We again examine $\langle \chi_s,\chi_s\rangle$ to see when we obtain $1$.
 Taking $g_t\in (a^t)^G$ an argument similar to the $\theta_j$ case gives   
 \begin{align*}
 &\langle \chi_s\uparrow G, \chi_s  \uparrow G\rangle
 =\frac{1}{|G|} \sum_{g\in G} \left(\chi_s\uparrow G\right)(g)\cdot\left(\chi_s\uparrow G\right)\left(g^{-1}\right)\\
 =&\frac{1}{2q^6-2q^2}\left(\sum_{g\in H} \left(\chi_s\uparrow G\right)^2(g)\right)\\
 =&\frac{1}{2q^6-2q^2}\left((2q^2+2)^2 + 4(q^4-1) + (q^4+q^2)\sum_{t=1}^{(q^2-2)/2}\left(\chi_s(g_t)\right)^2\right)\\
 =&\frac{1}{2q^6-2q^2}\left(8q^4+8q^2 + (q^4+q^2)\sum_{t=1}^{(q^2-2)/2}\left(\rho^{st} + \rho^{-st} + \rho^{stq} + \rho^{-stq}\right)^2\right)\\
 =&\frac{1}{2q^6-2q^2}\Bigg(8q^4+8q^2 + (q^4+q^2)\sum_{t=1}^{(q^2-2)/2}\Big(4+\rho^{2st} + \rho^{-2st} + \rho^{2stq} + \rho^{-2stq}\\
  &\qquad \qquad\qquad\qquad\qquad\qquad\qquad\qquad\qquad + 2\rho^{-st(q-1)} + 2\rho^{-st(q+1)}\Big)\Bigg)\\
 =& \frac{1}{2q^6-2q^2}\left(8q^4+8q^2 + \left(2(q^2-2)-2\right)(q^4+q^2) + (q^4+q^2)\sum_{t=1}^{q^2-2}(2\rho^{st(q+1)} + 2\rho^{st(q-1)})\right)\\
 =& \frac{1}{2q^6-2q^2}\left(8q^4+8q^2 + 2q^6-4q^4-6q^2 +  (q^4+q^2)\sum_{t=1}^{q^2-2}(2\rho^{st(q+1)} + 2\rho^{st(q-1)})\right)\\
 =& \frac{1}{2q^6-2q^2}\left(2q^6+4q^4+2q^2 +  (q^4+q^2)\sum_{t=1}^{q^2-2}(2\rho^{st(q+1)} + 2\rho^{st(q-1)})\right).
 \end{align*}

Here we used the facts that $(q^2-1)\nmid 2s$ and $(q^2-1) \nmid 2sq$, since $1\le s\le (q^2-2)/2$. Now, since $(q^2-1) \nmid s$, only one of $(q^2-1) \mid s(q+1)$ or $(q^2-1) \mid s(q-1)$ can be true, this shows that the above is equal to 
$$ \begin{cases}
\frac{1}{2q^6-2q^2}\left(2q^6+4q^4+2q^2 -4(q^4+q^2)\right) = \frac{2q^6-2q^2}{2q^6-2q^2} =1 \text{ if } (q^2-1) \nmid s(q\pm 1)\\
\frac{1}{2q^6-2q^2}\left(2q^6+4q^4+2q^2 +  2(q^4+q^2)(q^2-3)\right) = \frac{4q^6-4q^2}{2q^6-2q^2} = 2 \text{ if } (q^2-1) \mid s(q\pm 1).
\end{cases}
$$

Since there are $\frac{q}{2}$ values of $s$ for which $(q^2-1) \mid s(q+1)$ and $\frac{q-2}{2}$ values where $(q^2-1) \mid s(q-1)$,  we see  that $\frac{2q-2}{2}=q-1$ characters $\chi_s$ of $H$ split in $G$. Then the remaining $\frac{q^2-2q}{2}$ characters fuse in $G$. Recall that $1\leq j \leq q^2/2$ and $1 \leq s \leq (q^2-2)/2$. 
So 
$$
\deg\left(\tau_{G}\right) = 2 + 2q^2 + \left(2(q-1) + \frac{q^2-2q}{2}\right)(q^2+1) + \frac{q^2}{2}\left(q^2-1\right) = q^4+q^3+q.
$$

By Theorem  \ref {lem:totalsymplectic}  $q^4+2q^3+2q^2+2q+1$ is the largest degree of an irreducible character of $\mathrm{Sp}_4(q), q\geq 4$.   Since $G = \mathrm{Sp}_2\left(q^2\right)\colon 2$ and
$$
\deg\left(\tau_{G}\right) =  q^4+q^3+q < q^4+2q^3+2q^2+2q+1
$$
by Lemma \ref {lem:totalchar}  $\left(\mathrm{Sp}_4(q), \mathrm{Sp}_2\left(q^2\right)\colon 2\right)$ is not a strong Gelfand pair.
\qed
\medskip

Similar to Corollary \ref {cor1} we see that by Lemma \ref {isomorpicmaximal} and the fact that the above argument is a `total character argument' (not dependent on the particular embedding of $ \mathrm{Sp}_2(q^2): 2$ in $\mathrm{Sp}_4(q)$)        we have:

\begin{cor} \label  {Cor2} The maximal subgroup {\rm SO}$^-_4(q)<{\rm Sp}_4(q)$ is not a strong Gelfand subgroup.\qed
\end{cor}\medskip

\begin{thm}
For $q=2^e, e>1$, and $q_0$ such that $q=q_0^r$ for a prime $r$, the maximal subgroup $\mathrm{Sp}_4(q_0)\leq \mathrm{Sp}_4(q)$ is not a strong Gelfand subgroup.
\end{thm}
\noindent\textit{Proof.} 
By Theorem \ref {lem:totalsymplectic}  $\deg(\tau_{\mathrm{Sp}_4(q)}) = q^6+q^4-q^2$ for all even $q$. Then $\deg(\tau_{\mathrm{Sp}_4(q_0)}) = q_0^6+q_0^4-q_0^2$ and since $q=q_0^r$ and $r\geq 2$ we  see that
$$
q^4+q^3+q^2+q = q_0^{4r}+q_0^{3r}+q_0^{2r}+q_0^r \geq q_0^8+q_0^6+q_0^4+q_0^2.
$$
This shows  that
$$
\deg(\tau_{\mathrm{Sp}_4(q_0)}) = q_0^6+q_0^4-q_0^2 < q^4+2q^3+2q^2+2q+1
$$
and so by Lemma \ref{lem:totalchar} $(\mathrm{Sp}_4(q), \mathrm{Sp}_4(q_0))$ is not a strong Gelfand pair. \qed\medskip

\begin{thm}\label{thm:suz}
For $q=2^{2n+1}$, with $n$ a positive integer, the maximal subgroup $\mathrm{Sz}(q)$ in $\mathrm{Sp}_4(q)$ is not a strong Gelfand subgroup.
\end{thm}
\noindent\textit{Proof.} In \cite{Suz} Suzuki gives the irreducible characters of $\mathrm{Sz}(q)$, where $q=2^{2n+1}$. They are:

\noindent 
($i$)  the trivial character of degree $1$;

\noindent 
($ii$)  a doubly transitive character of degree $q^2$;

\noindent 
($iii$)  $(q-2)/2$ characters of degree $q^2+1$;

\noindent 
($iv$)  two complex characters of degree $2^n(q-1)$;

\noindent 
($v$)  $(q+2^{n+1})/4$ characters of degree $(q-2^{n-1}+1)(q-1)$;

\noindent 
($vi$)  $(q-2^{n+1})/4$ characters of degree $(q+2^{n-1}+1)(q-1)$.

This gives the following expression  for $\deg(\tau_{\mathrm{Sz}(q)})$:
\begin{align*}
&1+q^2+\left(\frac{q-2}{2}\right)\left(q^2+1\right) + 2\cdot 2^n\left(q-1\right)
+ \left(\frac{q+2\cdot 2^n}{4}\right)\left(q-2\cdot 2^n+1\right)\left(q-1\right)\\
&\qquad + \left(\frac{q-2\cdot 2^n}{4}\right)\left(q+2\cdot 2^n+1\right)\left(q-1\right)\\
& = 2^{n+1}(q-1)-q(q-1)+q^3.
\end{align*}
We now notice that the degree of the total character of $\mathrm{Sz}(q)$ is smaller than the maximal degree of an irreducible character in $\mathrm{Sp}_4(q)$ by Theorem  \ref{lem:totalsymplectic}. This shows by Lemma \ref{lem:totalchar} that $(\mathrm{Sp}_4(q), \mathrm{Sz}(q))$ is not a strong Gelfand pair when $q=2^{2n+1}$. \qed\medskip

This completes consideration of all the  maximal subgroups
listed in Table 1 and so concludes the proof of Theorem 
\ref {Thm1}.

\section*{Appendix}
\begin{lstlisting}
IsStrongGelfandPair := function(g, h);
	tf := true;
	ctg := CharacterTable(g);
	cth := CharacterTable(h);
	for character in ctg do
		r := Restriction(character, h);
		for i := 1 to #cth do
		       if InnerProduct(r, cth[i]) gt 1 then 
				tf :=false;
				break character;
		       end if;
		end for;
	end for;
	return tf;
end function;


G := SymplecticGroup(4,4);
[IsStrongGelfandPair(G, u`subgroup) : 
				u in MaximalSubgroups(G)];
\end{lstlisting}

\noindent
{\bf Acknowledgment} All computations made in the writing of this paper were accomplished  using Magma \cite{Mag}. Thanks are due to some anonymous referees for helpful comments.


\begin{thebibliography}{HJ}

\bibitem{A}
 Aizenbud, A.; Gourevitch, D.; Rallis, S.; Schiffmann, Gerard G., \emph{Multiplicity one theorems}, Ann. Math. (2) 172(2) (2010) 1407--1434.
\bibitem{AA} Aizenbud, A.; Gourevitch, D., \emph{Multiplicity one theorem for $({\rm GL}_{n+1} (\mathbb R), {\rm GL}_n (\mathbb R))$}, Selecta Math. (N.S.) 15(2) (2009) 271--294.

\bibitem{AHN} Anderson, Gradin; Humphries, Stephen P.; Nicholson, Nathan \emph{Strong Gelfand pairs of symmetric groups.}
J. Algebra Appl. 20 (2021), no. 4, Paper No. 2150054, 22 pp.

 \bibitem{BH} Barton, Andrea; Humphries, Stephen \emph{Strong Gelfand Pairs of {\rm SL}($2,p$)}. J. Algebra Appl. 22 (2023), no. 6, Paper No. 2350133, 13 pp.  https://doi.org/10.1142/S0219498823501335


\bibitem{max}
 Bray, John N.; Holt, Derek F.; Roney-Dougal, Colva M. \emph{The maximal subgroups of the low-dimensional finite classical groups. With a foreword by Martin Liebeck.} London Mathematical Society Lecture Note Series, 407. Cambridge University Press, Cambridge, 2013. MR3098485
\bibitem{BC} Brou, Kouakou Germain; Coulibaly, Pie; Kangni, Kinvi, 
\emph{Generalized Gabor transform for a strong Gelfand pair.} J. Adv. Math. Stud. 18 (2025), no. 1, 109--121.



\bibitem{Mag}  Bosma, Wieb; Cannon, John;    Playoust, Catherine; \emph{The Magma algebra system. I. The user language}, J. Symbolic Comput., 24 (1997), 235--265.
\bibitem{C} Chan, Kei Yuen,   \emph{Ext-multiplicity theorem for standard representations of $({\rm GL}_{n+1},{\rm GL}_{n})$.} Math. Z. 303 (2023), no. 2, Paper No. 45, 25 pp.

\bibitem{Har} Ceccherini-Silberstein, T.; Scarabotti, F.; Tolli, F., \emph{Harmonic analysis on finite groups, in Representation Theory, Gelfand Pairs and Markov Chains}, Vol. 108 Cambridge Studies in Advanced Mathematics, Cambridge University Press, Cambridge, 2008, xiv + 440 pp.

\bibitem{Atl} Conway, J. H.; Curtis, R. T.; Norton, S. P.; Parker, R. A.; Wilson, R. A. \emph{Atlas of finite groups. Maximal subgroups and ordinary characters for simple groups.} With computational assistance from J. G. Thackray. Oxford University Press, Eynsham, 1985. xxxiv+252 pp.

\bibitem{Da} 
Dabbaghian-Abdoly, Vahid  
\emph{Characters of some finite groups of Lie type with a restriction containing a linear character once}. 
J. Algebra 309 (2007), no. 2, 543--558.



\bibitem{Dor} Dornhoff, L (1971) \emph{Group Representation Theory. Part A: Ordinary Representation Theory}. Pure and Applied Mathematics. Vol. 7. New York: Marcel Dekker, Inc., pp. vii+254. [Google Scholar]

\bibitem{speven}
Enomoto, Hikoe \emph{The characters of the finite symplectic group} ${\rm Sp}(4,\,q)$, $q=2\sp{f}$. Osaka Math. J. \textbf{9} (1972), 75--94. MR0302750

 

\bibitem{GH}
Gardiner, Jordan C.; Humphries, Stephen P. \emph{Strong Gelfand pairs of ${\rm SL}(2,p^n)$}. Comm. Algebra 52 (2024), no. 8, 3269--3281. 


\bibitem{Tot2} Humphries, Stephen; Kennedy, Chelsea; Rode, Emma \emph{The total character of a finite group}. Algebra Colloq. 22 (2015), Special Issue no. 1, 775-778. (Reviewer: Silvio Dolfi)

 

\bibitem{JaL} 
James, Gordon, and Liebeck, Martin. \emph{Representations and Characters of Groups}. 2nd ed., Cambridge University Press, 2001.



\bibitem{JKar}Karlof, John \emph{The Subclass Algebra Associated with a Finite Group and Subgroup}, Amer. Math. Soc. Vol 207 (Jun., 1975) pp. 329-341.

  

\bibitem{Tot1} Prajapati, S. K.; Sarma, R. \emph{Total character of a group $G$ with $(G, Z(G))$ as a generalized Camina pair}. Canad. Math. Bull. 59 (2016), no. 2, 392--402.

\bibitem{Tot3} Prajapati, S. K.; Sury, B. \emph{On the total character of finite groups.} Int. J. Group Theory 3 (2014), no. 3, 47-67.

 


\bibitem{Suz}
Michio Suzuki. \emph{A new type of simple groups of finite order}. Proc. Nat. Acad. Sci. U.S.A.
46 (1960), 868--870. MR0120283
\bibitem{Den} Travis, Dennis \emph{Spherical Functions on Finite Groups}, Journal of Algebra 29, 65-76 (1974).


\bibitem{Wik} Wikipedia contributors. ``Gelfand pair." Wikipedia, The Free Encyclopedia, 25 May 2021.







\end{thebibliography}
\end{document}